\pdfoutput=1
\documentclass[12pt,a4paper,reqno]{amsart}
\usepackage{amssymb}

\usepackage{amscd}
\usepackage{enumerate}
\usepackage{graphicx}
\usepackage{siunitx}
\usepackage{tikz-cd}
\usepackage{color}
\usetikzlibrary{arrows}
\numberwithin{equation}{section}

\usepackage{mathabx}

\usepackage{mathtools}
\usepackage[tableposition=top]{caption}
\usepackage{booktabs,dcolumn}

\DeclareCaptionLabelSeparator{separation}{:\quad}
\DeclareFontFamily{OT1}{rsfs}{}
\DeclareFontShape{OT1}{rsfs}{n}{it}{<-> rsfs10}{}
\DeclareMathAlphabet{\mathscr}{OT1}{rsfs}{n}{it}

\theoremstyle{plain}

\theoremstyle{definition}

\renewcommand{\phi}{\varphi}

\usepackage{hyperref}

\usepackage{algorithm, algorithmic}

\begin{document}

\title{The ``bounded gaps between primes'' Polymath project - a retrospective}

\author{D.H.J. Polymath}
\address{http://michaelnielsen.org/polymath1/index.php}

\begin{abstract}  For any $m \geq 1$, let $H_m$ denote the quantity $H_m := \liminf_{n \to \infty} (p_{n+m}-p_n)$, where $p_n$ denotes the $n^{\operatorname{th}}$ prime; thus for instance the twin prime conjecture is equivalent to the assertion that $H_1$ is equal to two.  In a recent breakthrough paper of Zhang, a finite upper bound was obtained for the first time on $H_1$; more specifically, Zhang showed that $H_1 \leq \num{70000000}$.

Almost immediately after the appearance of Zhang's paper, improvements to the upper bound on $H_1$ were made.  In order to pool together these various efforts, a \emph{Polymath project} was formed to collectively examine all aspects of Zhang's arguments, and to optimize the resulting bound on $H_1$ as much as possible.  After several months of intensive activity, conducted online in blogs and wiki pages, the upper bound was improved to $H_1 \leq \num{4680}$.  As these results were being written up, a further breakthrough was introduced by Maynard, who found a simpler sieve-theoretic argument that gave the improved bound $H_1 \leq \num{600}$, and also showed for the first time that $H_m$ was finite for all $m$.  The polymath project, now with Maynard's assistance, then began work on improving these bounds, eventually obtaining the bound $H_1 \leq \num{246}$, as well as a number of additional results, both conditional and unconditional, on $H_m$.

In this article, we collect the perspectives of several of the participants to these Polymath projects, in order to form a case study of online collaborative mathematical activity, and to speculate on the suitability of such an online model for other mathematical research projects.
\end{abstract}


\maketitle

\section{Introduction}

In a blog post \cite{gowers-blog} from January 2007 entitled ``Is massively collaborative mathematics possible?'', Timothy Gowers proposed a new format for solving mathematical research problems, based on a large number of small contributions from many mathematicians in a public online forum, as opposed to the traditional model of a small number of mathematicians collaborating privately and intensively on a single problem.  After an enthusiastic online response to this blog post, Gowers soon launched the first of the ``Polymath'' projects.  This project, now known as Polymath1, had a specific, and feasible-looking goal: to take an important theorem in density Ramsey theory (the density Hales-Jewett theorem of Furstenberg-Katznelson \cite{furst}), which until this project only had a difficult ergodic-theoretic proof, and find a purely combinatorial proof of this statement.  Within 24 hours of the online launch of this project, Gowers' blog post was already flooding with a large number of mathematical comments, ideas, and observations.  Promising approaches were soon identified and explored, reading seminars to try to digest difficult papers in the literature were set up, and a wiki to maintain links to all relevant online resources, as well as lengthier arguments not easily placed in a blog comment, was also set up.  After six weeks of hectic activity, with contributions from dozens of mathematicians, Gowers was able to declare victory: a convincing outline of an elementary proof of the density Hales-Jewett theorem had been found.  The process of formally writing up the results \cite{Polymath1a,Polymath1b} took somewhat longer (closer to six months than six weeks); nevertheless, the project was a considerable success by most reasonable standards.  See \cite{nature} for a further discussion of the Polymath1 project.

After Polymath1, several further Polymath projects (and Polymath-like projects) were launched by Gowers and several other mathematicians, with varying degrees of success.  For instance, Polymath4, whose aim was to find fast deterministic methods to locate primes, ended up with a modest partial result, published\footnote{The original intention was to publish this article under the agreed-upon pseudonym of ``D.H.J. Polymath'', but the journal requested that the authors' true names and affiliations be used instead.} in \cite{Polymath4}, although the original aim of improving the fastest previously known such method was not actually achieved.  An online collaborative project to analyse a claimed proof of the infamous $P \neq NP$ problem by Deolilakar, was not ``officially'' a Polymath project, but followed its format (and used the same wiki) as the Polymath projects; it ended up with a strong consensus that the errors discovered in the argument were unlikely to be fixable; see \cite[Chapter 1]{lipton}.  A number of smaller ``mini-Polymath'' projects, based around solving a problem from the International Mathematical Olympiads, were also successful (with the problem chosen invariably being solved with 48 hours) and generally enjoyed by the participants involved, although such projects were more of a social mathematical event than a genuine research mathematical program; see \cite{pease} for an analysis of one of these projects.  On the other hand, several Polymath projects, after promising starts (with a level of intensity and participation comparable to that of the Polymath1 project), eventually stalled, with initial promising ideas to attack the problem encountering seemingly insuperable obstacles, and with the level of attention given by the participants falling below critical mass.  However, even the stalled projects have made some contributions towards their objectives, which (due to the open nature of these projects) are freely available for other researchers to use; for instance, some of the ideas developed in the stalled Polymath7 project to solve the ``hot spots'' conjecture for acute triangles were used (with permission) in a separate publication \cite{siudeja} outside of the Polymath project.

This article is devoted to the most recent Polymath projects, namely the Polymath8 project to understand, build upon, and improve the breakthrough work of Zhang \cite{zhang} on bounded gaps between primes. More precisely, for any $m \geq 1$, let $H_m$ denote the quantity $H_m := \liminf_{n \to \infty} (p_{n+m}-p_n)$, where $p_n$ denotes the $n^{\operatorname{th}}$ prime; thus for instance the twin prime conjecture is equivalent to the assertion that $H_1$ is equal to two.  Prior to the work of Zhang, it was not known unconditionally if any of the $H_m$ were finite, although under the strong additional assumption of the Elliott-Halberstam conjecture \cite{elliott}, the bound $H_1 \leq 16$ was shown in an important paper of Goldston, Pintz, and Y{\i}ld{\i}r{\i}m \cite{gpy}, who also established the unconditional result $\liminf_{n \to \infty} \frac{p_{n+1}-p_n}{\log p_n} = 0$.  By reviving and brilliantly modifying some arguments of Bombieri, Fouvry, Friedlander, and Iwaniec \cite{fi, fi-2, bfi, bfi-2, bfi-3}, and combining\footnote{For a more in-depth treatment of the mathematics of Zhang's work, see \cite{granville}.} this with variants of the Goldston-Pintz-Y{\i}ld{\i}r{\i}m argument, Zhang established the bound
$H_1 \leq {70000000}$.  In his paper, Zhang comments: ``This result is not optimal ... to replace the [upper bound on $H_1$] by a value as small as possible is an open problem that will not be discussed in this paper''.

This provided an irresistible challenge to many of the mathematicians following these developments.  One week after the release of Zhang's preprint, it was remarked by Mark Lewko (as an offhand comment on the web site MathOverflow) that by replacing a single crude upper bound in Zhang's paper, one could improve the upper bound on $H_1$ to ${63,374,611}$; one week after that, it was observed by Trudgian \cite{trudgian}, in a short arXiv preprint, that another ``cheap'' tweak to Zhang's argument could lower the bound further to ${59,874,594}$.  A few days later, in a blog post of Morrison \cite{morrison}, another slight (and computer-assisted) optimization to Zhang's construction improved the bound further to ${59,470,640}$.  At this point, several other mathematicians began to join in the comments of Morrison's post, contributing further small observations and improvements - very much in the spirit of the ``Polymath'' enterprise.  Indeed, after some quick discussion on the Polymath blog \cite{Polymath}, it was agreed to retroactively convert Morrison's post to the first post on what would now become the eighth Polymath project, whose stated aims were to understand Zhang's argument as thoroughly as possible, and then to optimize the bound on $H_1$ from that argument.  After three months of intense activity, in which the upper bound $H_1$ was repeatedly lowered\footnote{A timeline of improvements may be found at \url{michaelnielsen.org/polymath1/index.php?Timeline_of_prime_gap_bounds}.} by a variety of methods, the upper bound on $H_1$ eventually stabilized at $H_1 \leq {4680}$.  

At this point, attention began to turn to the lengthy, but necessary, process of writing up the results \cite{Polymath8a}.  But as the paper was nearing completion, news of another breakthrough in the subject emerged: by a clever modification to the original sieve-theoretic techniques of Goldston, Pintz, and Y{\i}ld{\i}r{\i}m, Maynard \cite{maynard} was able to circumvent many of the technicalities of the arguments of Zhang (or of Polymath8), and obtain a further significant improvement of the bound, to $H_1 \leq {600}$; furthermore, bounds on $H_m$ for larger $m$ were now available, and under the assumption of the Elliott-Halberstam conjecture, the previous bound of $H_1 \leq 16$ was improved slightly to $H_1 \leq 12$.

With this new advance, the Polymath8 participants decided to start a second phase to the project (concurrently with finishing the writing up of the first phase).  In this second phase, dubbed Polymath8b, the new methods of Maynard would be combined with the existing results of the first phase (now dubbed Polymath8a) to obtain further improvements to $H_1$ and $H_m$.  After about eight months of intense activity, the results of this paper were also written up and submitted \cite{Polymath8b}; the bound on $H_1$ had been lowered unconditionally to $246$, and the bound $H_1 \leq 6$ obtained assuming the generalized Elliott-Halberstam conjecture (with the latter bound being the limit of the sieve-theoretic method).  Furthermore, a number of technical sieve-theoretic innovations have been discovered en route to these results, which will hopefully have value for other problems in analytic number theory beyond the bounded gaps between primes problem.

In this retrospective, we record various first-hand impressions of the Polymath8 project as it unfolded, both from active participants and from casual followers, who have reported both on their experience and on their thoughts on whether the success of this particular project can be replicated for other problems.  The contributions are arranged in chronological order of submission.

\section{Terence Tao}

\emph{Terence Tao is a professor of mathematics at UCLA.}

The first that I heard about Zhang's breakthrough on bounded gaps between primes was on May 13, 2013, when I started getting emails from friends and colleagues about a rumoured breakthrough on this problem in number theory.  I was intrigued, but initially sceptical; the closely related twin prime conjecture, for instance, has attracted more than its fair share of arguments that ended up completely falling apart, and Yitang Zhang had published little in the area for over a decade.  But then a trusted friend passed on the news that Zhang's result on bounded gaps had apparently already been refereed and accepted by the Annals of Mathematics, and also reported some details of the argument (a combination of the Goldston-Pintz-Y{\i}ld{\i}r{\i}m method and the arguments of Bombieri, Fouvry, Friedlander, and Iwaniec).  At that point, I felt the information I had was credible enough to post a brief comment about it on my own Google+ page.  By the next day, the Annals put up the preprint on its web page, Zhang presented his work in a special seminar in Harvard, and suddenly the online maths community was buzzing with the result.  

Initially, I was preoccupied with other research projects, and felt content to let the rest of the analytic number theory community digest the result; I sent Zhang a quick congratulatory email and returned to work.  A week or two later, when the first blog post appeared in which Zhang's original bound of $\num{70000000}$ was lowered to $\num{59470640}$ with some computer assistance, I again commented on this on my Google+ page, figuring that this could be an opportunity for number theory enthusiasts to make a nice contribution to this subject, and planned to again return to work... but then decided to spend just a few minutes fiddling with the problem myself.  

Zhang arrived at the bound of $\num{70000000}$ by constructing an ``admissible $k$-tuple'' whose diameter was bounded by this amount.  An admissible $k$-tuple is a set of $k$ increasing integers, with the property that it misses at least one residue class mod $p$ for each prime $p$; for instance $(0,2,6)$ is an admissible $3$-tuple, but $(0,2,4)$ was not.  These tuples play a central role in the famous \emph{prime tuples conjecture} of Hardy and Littlewood \cite{hardy}, which generalizes the twin prime conjecture and is still completely open, despite many partial results.  In Zhang's work, $k$ had to be sufficiently large; his analysis allowed him to take $k = \num{3500000}$.  Zhang constructed his admissible tuple by considering a block of $k$ consecutive primes that was guaranteed to be admissible; the subsequent improvements came from finding slightly narrower blocks of consecutive primes and checking the admissibility by computer.  The discussion reminded me of some work of Hensley and Richards \cite{hensley} who used admissible tuples to show that the prime tuples conjecture contradicted another conjecture of Hardy and Littlewood about the prime counting function $\pi(x)$; I did a quick calculation that suggested that this result ought to lower the bound further to $\num{58885998}$, posted this on the blog, and then tried to return to ``real'' work again.

But something kept drawing me back.  At the time, I was working on a lengthy research project which I felt would take months, if not years, to complete (and indeed, it is still far from completion as of this time of writing, and in fact was shelved as the Polymath project began to absorb more of my time).  In contrast, with just ten minutes of effort, there was a chance to push the bound down a little more and claim, however briefly, the ``world record'' for bounded gaps betweeen primes.  Plus, it gave me an excuse to actually go through Zhang's paper in more detail.  So I found myself returning again to the blog (which was now becoming quite lively with other comments and contributions), and digging further and further into Zhang's paper to try to squeeze out more improvements.  The initial gains had come purely from reading the first two pages of Zhang's 56-page paper, keeping his value $\num{3500000}$ of $k$ intact; by reading up the next twenty or so pages, we found out how this number was arrived at (by combining a sieve-theoretic argument of Goldston, Pintz, and Y{\i}ld{\i}r{\i}m with a deep new distributional estimate on primes), and began to optimise this value too.  Again, we began with ``cheap'' numerical optimisations, for instance using numerics to lower $k$ to $\num{2947442}$, which in turn lowered the $H_1$ bound to $\num{48112378}$, but then started looking carefully at how Zhang controlled various error terms to make more significant gains; within a few days, we had reduced the bounds for both $k$ and $H_1$ by approximately a full order of magnitude.

By this point it had become abundantly clear that the project could benefit from being organised as a Polymath project, with a wiki page to keep track of all the progress and links to resources, and with a separate thread to discuss administrative issues.  We hastily proposed converting Morrison's blog post retroactively to such a project, with the twin goals of improving the bound on $H_1$, and also to understand and clarify Zhang's argument (and related arguments in the literature).  Reaction to this proposal was generally positive, but there were some objections.  Chief among these was the concern that the existence of such a prominent project, with many mathematicians participating, might discourage and intimidate others from trying to contribute to this development of the subject.  However, it was felt that there would be a large surge of interest in the problem even if there was no formal Polymath project, and that the narrow focus of the project would allow other work relating to Zhang's breakthrough.  Indeed, in the months after Zhang's paper, there were a number of other preprints \cite{pintz-polignac}, \cite{anders}, \cite{maynard-ap}, \cite{freiberg}, \cite{banks}, \cite{thorner} that improved Zhang's result in other ways than a direct improvement of $H_1$, as well as preprints \cite{pintz-improv}, \cite{maynard} that improved the $H_1$ bound beyond the record obtained by Polymath8; so it appears that the Polymath8 project did not, in fact, ``crowd out'' traditional research in this direction.  Indeed, with the additional breakthrough of Maynard \cite{maynard}, the pace of new research papers in this area has even accelerated \cite{consecutive}, \cite{benatar}, \cite{li-pan}, \cite{castillo}, \cite{pollack}, \cite{banks2}, \cite{clusters}, \cite{lola}, \cite{pintz-ratio}, \cite{chua}, \cite{pintz-new}. 

The Polymath8 project naturally split into three pieces: the search for narrow admissible tuples of a given cardinality $k$; the refinement of the value of $k$; and the further reading of Zhang's paper.  Scott Morrison's blog post was already hosting a lively discussion of the first component of this project; in early June, I wrote two further blog posts to host the discussion of the other two components, while Scott ``rolled over'' his blog post (as per the usual Polymath custom) to a fresh one, in which he summarised the previous progress.  

The most immediate route to improving $k$ was to digest and then optimise the sieve-theoretic components of Zhang's work; this followed fairly closely the previous work of Goldston, Pintz, and Y{\i}ld{\i}r{\i}m which I was already familiar with, so I spent most of my efforts on this part of the project, in particular writing up detailed notes on the GPY sieve on my blog.  Almost immediately, high-quality mathematical comments started rolling in from many different readers: typos and other minor errors in the blog post were quickly pointed out, and different aspects of the argument were discussed.  It was quickly noticed that a certain cutoff function $f$ in the definition of the sieve could be optimised to improve the value of $k$; this led to a calculus of variations problem which was soon solved with the aid of Bessel functions, bringing $k$ down (in principle, at least) all the way down to $\num{34429}$ (and $H_1$ to about $\num{390000}$).  It was then quickly pointed out that this Bessel function optimisation had been previously worked out in \cite{farkas} (and in earlier unpublished work of Conrey); indeed, one of the advantages of the Polymath projects, with their broad level of participation, is that connections to relevant literature are very likely to be unearthed by at least one of the participants.

By the second week of June, a couple of us had managed to complete our collaborative online reading of Zhang's paper, in particular understanding how the distributional estimates on primes that were used to bound $k$ arise from three key estimates, that Zhang called ``Type I'', ``Type II'', and ``Type III'' estimates.  I then prepared some technical blog posts on these estimates, as well as the combinatorial argument Zhang used to merge them together.  The Polymath project had now split into \emph{five} active and loosely interacting components: one group was focusing on optimising the Type I and Type II estimates (which were a lengthy but fairly straightforward application of standard analytic number theory methods, such as the Linnik dispersion method and Weil exponential sum estimates); one group on the Type III estimates (which were somewhat more exotic, relying on an exponential sum estimate of Birch and Bombieri which in turn used Deligne's work on the Riemann hypothesis on varieties); one group on the combinatorics of putting these estimates together to create a distributional estimate on primes; one group on optimising the GPY sieve to convert distributional estimates on primes into a concrete value of $k$; and one group to find narrow prime tuples which converted values of $k$ to values of $H$.  We had managed to organise ourselves into a sort of factory production line: an advance in, say, the Type I estimates would be handed over to the combinatorics group to produce a new distributional estimate in primes, which the sieve team would then promptly convert into a revised value of $k$, which the prime tuples team would then use to update their value of $H_1$.  Steady improvements from all of these groups led to new improvements to these values on a daily basis for several weeks (see \url{michaelnielsen.org/polymath1/index.php?title=Timeline_of_prime_gap_bounds} for a timeline of these gains).  

Many important contributions came from mathematicians who were not initially involved in the Polymath project.  For instance, in the first week of June, Janos Pintz released a preprint \cite{pintz-improv} in which he went through Zhang's paper and optimised all the numerical constants, leading to further improvements in the values of $k$ and $H_1$; within a day or two, the Polymath team had gone through the preprint (correcting some arithmetic errors along the way), lowering $k$ to $\num{26024}$ and $H_1$ to about $\num{280000}$.   Later, we ran into a problem in which we became victims of our own success: the values of $k$ had dropped so low (around $\num{6000}$) that a certain error term in the sieve analysis that decayed exponentially in $k$ became non-negligible.  Pintz, who had followed developments closely, came up with a much more efficient way to truncate the sieve that rendered this error term almost non-existent (and dropping $k$ almost immediately from $\num{6000}$ to about $\num{5000}$), and then shared these computations with the Polymath project.   Similarly, Etienne Fouvry, Emmanuel Kowalski, Phillipe Michel, and Paul Nelson, who had previously worked on using Deligne-type exponential sum estimates to obtain distributional estimates similar to Zhang's Type III estimates, found a simpler approach than Zhang to these estimates that gave significantly better numerology, and they also donated these arguments to the Polymath project and began actively participating in the following discussion.  Thomas Engelsma, who had made extensive numerical computations of narrow admissible tuples over a decade ago, opened up his database to our project also.

By the end of June (with $k$ dipping below $\num{1000}$ and $H$ falling below $\num{7000}$), our understanding of various components of the project had matured significantly, though the pace of discussion was still lively (with perhaps twenty-odd comments each day, down from a peak of fifty or so per day).  A database had been set up at \url{math.mit.edu/~primegaps/} to automatically record the narrowest known prime tuples for a given value of $k$ (up to $k=5000$), which automated the task of converting $k$-values to $H$-values, with several of the more computationaly minded participants competing to add ever narrower tuples to the database.  A prime distribution estimate could similarly be converted into a value of $k$ with a few lines of Maple code, using Pintz's version of the truncated sieve, and a further few lines of code allowed us to automatically convert any improvement in the Type I, II, III estimates to a prime distribution estimate in an optimised fashion.  The Type II and Type III estimates had been optimised to such an extent that they were no longer the dominant barrier in further improvement of the prime distribution estimates.  The only remaining major source of improvement was the Type I estimates, but over the course of the next few weeks, several new ideas on optimising these estimates (most notably, using the $q$-van der Corput method of Graham and Ringrose \cite{graham}, combined with the theory of trace weights arising from the work of Deligne, and also carefully splitting the summations to maximise the effectiveness of the Cauchy-Schwarz inequality) were implemented, leading to the final values of $k=632$ and $H_1=4680$.

Further numerical improvements looked very hard to come by, and so for the next few months, the Polymath project devoted their efforts to writing up the results, splitting the paper into several sections on a shared Dropbox folder, and with various participants working in separate sections in parallel and coordinating their efforts through blog comments.  Somewhat to our surprise, once all the various arguments spread out over over a dozen blog posts had been collated together, the final paper turned out to be extremely large in size -- $163$ pages, to be exact\footnote{After significant revisions in the light of the referee report and subsequent developments, the paper has since been shortened somewhat to $107$ pages.} (for comparison, Zhang's original paper was $53$ pages in length).  Somehow, the Polymath format had managed to efficiently segment the research project into more manageable chunks, so that no one participant had to absorb the entire $163$-page argument at any given time while the research was ongoing.  

In late October, when the Polymath paper was going through several rounds of proofreading in preparation for submission, James Maynard announced a major breakthrough \cite{maynard}, in which the prime gap $H_1$ was lowered substantially, to $700$ (quickly revised downward to $600$), with the improvement $H_1 \leq 12$ available under the assumption of the Elliott-Halberstam conjecture.  Furthermore, Maynard had found a way to modify the GPY sieve argument in such a way that the new distribution estimates on primes (the most difficult and novel part of Zhang's work, and of the Polymath project also) were no longer needed, and furthermore the methods could for the first time be used to bound $H_m$ for larger $m$.  (I had also begun to arrive at similar conclusions independently of Maynard, though I did not get the precise numerical results that Maynard did.)  After some discussion with Maynard and with the Polymath team, it was decided that we should both publish our results separately, and then to join forces to see if the new ideas in Maynard's paper could be combined with those in the Polymath paper to obtain further improvements (this effort was soon dubbed ``Polymath8b'').

I had some mixed feelings about the continuation of the Polymath8 project; running the Polymath8a project had already consumed several months of my time, and I was thinking of turning to unrelated research projects.  However, the enthusiasm of the other participants and the lure of getting a quick payoff from comparatively brief snatches of mathematical thought were once again too great to resist.

Interestingly, Polymath8b ended up moving in a rather different mathematical direction than Polymath8a; it turned out to be quite difficult to effectively combine the prime distribution estimates from Polymath8a with Maynard's sieve to bound $H_1$ (although they could be used to bound $H_m$ for larger values of $m$), and it was soon realised that the most promising way forward was to either optimise or generalise a certain multidimensional calculus of variations problem used in Maynard's work.  

The former approach was initially more fruitful; by taking the code Maynard used to obtain lower bounds for a variational problem by testing that problem on polynomials, and making it run faster and more efficiently on more powerful computers, we managed to cut down the unconditional $H_1$ bound from $600$ to $300$, but improving the conditional bound of $H_1 \leq 12$ was far tougher; indeed, we had established an upper bound on the variational problem that Maynard used to show that we could not hope to improve upon this bound without modifying the sieve.  In order to break the $12$ barrier, the efficiency of our sieve (measured by a quantity that we called $M_4$, which arose from a four-dimensional variational problem) had to exceed $2$.  There was then a lengthy and frustrating ``Zeno's paradox'' period in which the efficiency $M_4$ of our sieves kept improving incrementally (from $1.845$, to $1.937$, to $1.951$, ...), but never quite enough to surpass the magic threshold of $2$ needed to break the barrier.  Finally, there was a breakthrough; after several weeks of effort, we stumbled upon an ``epsilon trick'' that allowed one to slightly enlarge the class of permissible cutoff functions in the variational problem, at the expense of worsening the quantity that one was trying to optimise.  It turned out that this tradeoff was advantageous, allowing us to move $M_4$ to $2.018$ and to reduce the bound on $H_1$ to $H_1 \leq 8$ (though we had to replace the Elliott-Halberstam conjecture with a strengthening of that conjecture, which we call the generalised Elliott-Halberstam conjecture).  A similar ``Zeno's paradox'' game then played out for the analogous three-dimensional variational problem $M_3$, which after several further refinements both to the sieve, and to the numerical procedure for optimising the sieve, eventually pushed $M_3$ to be larger than $2$ as well, giving the optimal result $H_1 \leq 6$ assuming the generalised Elliott-Halberstam conjecture.  (The parity obstruction of Selberg prevents any better bound on $H_1$ from purely sieve-theoretic considerations.)  Some of this new technology also allowed for some slight lowering of the unconditional bound of $H_1$ to $246$, but further improvement beyond this point seemed to require enormous amounts of computation, and by early May we were happy to ``declare victory'' at this point and write up the Polymath8b results.

The Polymath8 project had perhaps a dozen or two active participants, but many more mathematicians and interested amateurs followed the progress online.  During the project, whenever I visited another institution, I was usually asked what the latest value of $H_1$ was, and how low I thought it would go.  In that respect, we were fortunate that we had such a simple and easily understood statistic that could be used as a proxy for the degree of technical advance that we were making; it is not clear if future Polymath projects would be as easy to follow on a casual basis.

The project also forced me to think and work in different ways from what I was accustomed to.  I do not have extensive experience with programming, and most of the really heavy computational work was done by others on this project, but I did manage to write up some simple Maple code to at least verify the numerical computations that others had generated.  At another juncture, the way forward hinged on finding the optimal way to decompose the unit cube into polyhedral pieces; lacking sufficient geometric or algebraic intuition, I ended up having to build a cubic lattice out of my son's construction toys in order to visualise all the decompositions being proposed.  I also found myself having to learn areas of mathematics I would not otherwise have been exposed to, from $\ell$-adic cohomology to the Krylov subspace method.  All in all, it was an exhausting and unpredictable experience, but also a highly thrilling and rewarding one.



\section{Andrew Gibson}

\emph{Andrew Gibson is an undergraduate mathematics student at the University of Memphis.}

Shortly after Zhang announced his result and you (Tao) proposed the project, my classmates and I began a small weekly seminar with a professor devoted to studying some of the theory involved (analytic number theory, sieve methods, etc.), albeit on a much more elementary level that was within our reach. Of course, the majority of the actual proof is still mostly over our heads, but at least I feel as if I've gained a birds'-eye-view of the strategy and, probably more importantly, how it fits into the larger field. (For instance, before any of this, I could never have explained the Bombieri-Vinogradov theorem or the Hardy-Littlewood prime tuple conjecture.) So for us the project was a great excuse to enter a new subject, and has been immensely educational.

More than that, though, reading the posts and following the `leader-board' felt a lot like an academic spectator sport. It was surreal, a bit like watching a piece of history as it occurred. It made the mathematics feel much more alive and social, rather than just coming from a textbook. I don't think us undergrads often get the chance to peak behind closed doors and watch professional mathematicians ``in the wild'' like this, so from a career standpoint, it was illuminating. I get the sense that this is the sort of activity I can look forward to in grad school and as a post-doc doing research (...hopefully).

I also suspect that many other students from many other schools have had similar experiences but, like me, chose to stay quiet, as we had nothing to contribute. So, thank you all for organizing this project and for making it publicly available online.


\section{Pace Nielsen}

\emph{Pace Nielsen is an Assistant Professor at Brigham Young University.}

As a pre-tenure professional mathematician, I believe that a short introduction to my experiences in the Polymath8 project will be helpful to other young mathematicians deciding whether or not to participate in similar endeavors.  I initially joined the project for a number of reasons.  One is that I enjoy optimizing numerical results.  There is a certain pleasure that comes from deriving an elegant computational result.  Another reason I joined is that I wanted to understand Yitang Zhang's breakthrough as much as possible.  His result was quickly followed by James Maynard's (and independently Terry Tao's) multi-dimensional probabilistically motivated sieve, which is also incredibly interesting to me.  In my opinion, this improved sieve deserves a spot as one of the top advancements in analytic number theory of the past half-century.

My initial incursions into the project consisted mainly of pointing out minor corrections to the Polymath8a paper and asking questions about some of Tao's blog posts.  Hence, I don't consider myself a full participant in the 8a portion of the work; it was only during the 8b half of the work that I became a contributor.

There were a few things that surprised me about the whole experience. First was the friendliness of the other participants, particularly our host Terry Tao.  I want to publicly thank everyone involved in the project for the positive experience. A special thanks goes to James Maynard, whose kindness in sending me some of his original Mathematica code was what finally pushed me into full activity in the project.

Second, a large number of mathematicians I know commented (in personal communications to me) on the fact that they were ``impressed with my bravery'' in participating.  It caught me off guard that so many people had been following the project online, and that all of my comments (including my mistakes) were open to such a wide readership.  I believe it is important to consider this issue before deciding to participate in a public project. Some of the mistakes I made would never have seen the light of day in a standard mathematical partnership.  However, any collaboration relies on the ability for the participants to share ideas freely, even the ``dumb'' ones.

The third surprise was how much time I devoted to this project. I think that came from really getting excited about the work.  This is also the point at which I want to mention a (mild) concern on my part. As a pre-tenure professor I made the conscious choice going into this project that any time I spent was not necessarily going to be reflected in my tenure file. I knew that even if I contributed enough work to be noted as an official ``participant'' in the online acknowledgements, it wasn't clear whether this would count as co-authorship in the eyes of those on my tenure committee, or even of my colleagues generally.  Indeed, there are aspects of this type of cooperative mathematics which make it difficult to decide how participation in this venture should be treated by the mathematical community at large.  There are so many levels of effort that it can be somewhat confusing where the line is drawn between providing comments vs.\ being a full co-author vs.\ everything in-between. Also, the Polymath projects don't follow the convention used in other sciences of having the project manager decide who is a co-author on the final paper.

On the one hand, I support the idea of collaborative mathematics without an eye towards recognition.  With respect to the Polymath8b project in particular, I'm happy to give D.H.J.\ Polymath all the credit.  On the other hand, I do consider myself a co-author on the Polymath8b portion of the project, and thus want to take ownership for my part of the work.

The fourth and final surprise was that I did contribute something meaningful!  Sometimes this contribution happened by making simple comments on the blog.  For instance, I remember waking up one morning with the realization that a construction we were attempting would contradict the parity barrier in sieves.  This idea then led the experts to write some interesting formal mathematics on this issue.  Sometimes my contribution was made in the time-consuming busy-work of writing, running, and debugging computer code. And sometimes it was just in contributing my experience after thinking about the geometric picture long enough.  While it is intimidating working with Fields medalists and other experts, it is also a once in a lifetime opportunity to rub shoulders with such a wide array of mathematicians.

While I rate my involvement as extremely positive, others who are contemplating joining a massive mathematical online collaboration should keep in mind the costs in time, public embarrassment, and potential lack of control of your work before fully committing to the experience.


\section{James Maynard}

\emph{James Maynard is a Fellow at Magdalen College, Oxford.}

As a graduate student who had been looking at closely related problems, it was thrilling to hear of Zhang's initial breakthrough. I didn't participate in the subsequent Polymath 8a project, although I found myself reading several of the posts whilst studying Zhang's work for myself. I intended to avoid working on anything too close to the Polymath project (to avoid any competition), but one day I was going back through some ideas I'd had a several months earlier on modifications of the `GPY sieve', and I realized I could overcome the obstacles I had in my original attempt. This modification (also discovered by Tao) gave an alternative stronger approach to gaps between primes, although it didn’t produce the equidistribution results which lie at the heart of Zhang's work. With some small numerical calculations, this allowed me to show that there were infinitely many pairs of primes which differ by at most 600, and allowed one to show the existence of many primes in bounded length intervals.

I knew the numerical bounds in my work hadn't been fully optimized -- there was some slack in my approach, and there were also several opportunities to extend the method (such as incorporating ideas from Zhang and Polymath8a, or using more careful arguments). I was therefore pleased and excited when I learnt there was the intention for a Polymath8b project which I could be part of -- it was very exciting for me that there was such an interest in my work!

The style of a large collaborative project was very new to me. I had relatively little experience of research collaboration, and I was used to mainly trying out ideas alone. I certainly hadn't fully appreciated quite how public the posts were (and how many mathematicians who were not active participants would read the comments). In some ways this was quite fortunate; not realizing the attention posts might receive made me more willing to contribute openly. I posted several ideas which were not fully thought through, some of which were rather stupid in hindsight, but some of which I believe were useful to the project (and probably more useful than posting a fully thought out idea a few days later). The atmosphere of the project certainly helped encourage such partial contributions, which I feel was a large factor in the project's success.

It was remarkable (to me, at least) how smoothly the project went; this was partly to do with the problem having very clearly defined goals and being modular in its nature, but also due to the openness and friendliness of the participants. There seemed to be a good balance amongst the participants – some had more computational expertise, whilst others had a more theoretical background, and it was certainly useful to have both groups together for the effectiveness of the project. Much if the improvement in the unconditional bound came from extending the computations I had done initially, and it was certainly useful to have people who were rather more experienced than me with the larger computations.

I was surprised at how much time I ended up devoting to the Polymath project. This was partly because the nature of the project was so compelling – there were clear numerical metrics of `progress', and always several possible ways of obtaining small improvements, which was continually encouraging. The general enthusiasm amongst the participants (and others outside of the project) also encouraged me to get more and more involved in the project. Finally, the nature of the work also made it very suitable for working on in short bursts, which turned out to be very useful since I was travelling quite a lot whilst most of the project was underway.

I was aware that as a junior academic without a permanent position that I might ultimately not receive much credit in the eyes of hiring committees for participation in an atypical project such as the Polymath, where it is difficult to gauge the merit of my contribution. In my case, the fact that the project was so closely associated with my earlier work, and the fact that I found the project so interesting made me happy to accept this (although I made a conscious effort to continue to work on other projects at the same time). This is certainly something I feel any similarly junior prospective participant should be aware of, however.

Overall, I really enjoyed the Polymath experience. It was a great opportunity to work with several other mathematicians, and I feel pleased with the final results and my contribution to them.


\section{Gergely Harcos}

\emph{Gergely Harcos is a research advisor at Alfr\'ed R\'enyi Institute of Mathematics, and a professor at Central European University.}

I guess I am no longer a junior mathematician, which is a bad thing, but on the good side I can perhaps add a different perspective on my participation in the Polymath8 project.

Last year I applied for some serious grants, and it was on the same day when I learned that my proposals would be rejected. This was quite discouraging, and I felt as I needed to take some break from my main line of research. Around the same time, Zhang's exciting paper came out, and shortly after Terry Tao initiated a public reading seminar on his blog that later turned into the Polymath8 project. I was already familiar with the earlier breakthrough by Goldston-Pintz-Y{\i}ld{\i}r{\i}m; in fact I had incorporated it in my courses at Central European University and advised some students in related topics. I have always found this part of number theory very beautiful, although my research interests have been elsewhere. Following Terry's clean and insightful blog entries and the accompanying comments served for me two different purposes initially. First, I hoped to understand the new developments in a field that I found appealing. Second, I hoped to get a change of air in mathematics for the reasons explained above.

The Polymath8 project developed at blazing speed, and my initial goal was simply to catch up and read everything posted on the blog. This was quite challenging, because I am rather slow and prefer to check every line carefully, but at least I could serve as an early referee for the project. As a bonus, I got some ideas how to improve certain points in an argument already posted. In short, the Polymath8 project helped me to get going and feel myself useful, and participating was a lot of fun. At one point I embarrassed myself by posting several different ``proofs'' to an improved inequality that I conjectured, only to find out later that the claim was false. All this is recorded and preserved in the blog, but I do not regret it as it was honest and reflects the way mathematics is done. We try and we often fail.

It was also very interesting how my colleagues reacted. Some thought that one should not devote too much effort to a paper published under a pseudonym, but in fact my participation here got far more attention than elsewhere. For a couple of months the first question I was asked was about the current record on prime gaps. Another benefit of a Polymath project is that there is no pressure on the participants, one is free to join for a while then leave, and ignorance is normal as in a mathematical conversation.


\section{David Roberts}

\emph{David Roberts is a postdoctoral fellow at the University of Adelaide.}

I saw the announcement of Zhang's talk on Peter Woit's blog, and posted on Google+ (13 May 2013) about the twin prime problem, prime gaps more generally and about Zhang's talk and how big a deal it was. This post received much attention (more than I would have expected) and over the course of Polymath8, and the pre-official Polymath work, I kept posting the current records online, with explanation of what the progress meant or how it happened.

I'm not an analyst or number theorist (I work in category theory), so I was content to read the progress of the project and learn how all this business worked. I'd read through Zhang's preprint and was totally nonplussed, but the careful analysis -- and exposition! -- of Terry Tao and the other active participants made the ideas much clearer. In particular, concepts and tools that are well-known to analytic number theorists and are used without comment were brought into the open and discussed and explained.

When more specialised experts started working on subproblems, particularly numerical optimisation, it gave me snippets I could mention in my first-year algebra class to let them know how generalisations of the things they were learning (eigenvectors, symmetric matrices, convex optimisation etc) were being applied at the cutting edge of research. I even showed, on projector screens, Terry's blog and the relevant comments, some made that very day. It has been a great opportunity to expose students of all stripes to the idea of research in pure mathematics, and that a problem in number theory needed serious tools from seemingly unrelated areas.

For me personally it felt like being able to sneak into the garage and watch a high-performance engine being built up from scratch; something I could never do, but could appreciate the end result, and admire the process.


\section{Andrew Sutherland}

\emph{Andrew Sutherland is a Principal Research Scientist at MIT.}

I first heard about Zhang's result shortly after he spoke at Harvard in May, 2013. The techniques he used were well outside my main area of expertise, and I initially followed developments purely as a casual observer. It was only after reading Scott Morrison's blog, where people had begun discussing improvements to Zhang's bound, that I realized there was an interesting and essentially self-contained sub-problem (finding narrow admissible tuples) that looked amenable to number-theoretic combinatorial optimization algorithms, a subject with which I have some experience. I ran a few computations, and once I saw the results it was impossible to resist the urge to post them and join the Polymath8 project. Aside from interest in the prime gaps problem, I was curious about the Polymath phenomenon, and this seemed like a great opportunity to learn about it.

Like others, I was surprised by how much time I ended up devoting to the project. The initially furious pace of improvements and the public nature of the project made a very addictive combination, and I wound up spending most of that summer working on it. This meant delaying other work, but my collaborators on other projects were very supportive. I certainly do not begrudge the time I devoted to the Polymath8 effort; it was a unique opportunity, and I'm glad I participated.

In terms of the Polymath experience, there are a couple of things that stand out in my mind. In order to make the kind of rapid progress that can be achieved in a large scale collaboration, the participants really have to be comfortable with making mistakes in a forum that is both public and permanent. This can be a difficult adjustment, and there was certainly more than one occasion when I really wished I could retract something I had written that was obviously wrong. But one eventually gets used to working this way; the fact that everybody else is in the same boat helps. Actually, I think being forced to become more comfortable with making mistakes can be a very positive thing. This is how we learn.

The other thing that impressed me about the project is the wide range of people that made meaningful contributions. Not only were there plenty of participants who, like me, are not experts in analytic number theory, there were at least a few for whom mathematics is not their primary field of research. I think this is major strength of the Polymath approach, it facilitates collaboration that would otherwise be very unlikely to occur.

It is perhaps worth highlighting some of the features of this project that made it a particularly good Polymath candidate. First, the problem we were working on was well known and naturally attracted a lot of interested observers; this made it easy to recruit participants. Second, we had a clearly defined goal (improving the bound on prime gaps), and a metric against which incremental progress could be easily measured; this kept the project moving forward with lot of momentum. Third, the problem we were working on naturally split into sub-problems that were more or less independent; this allowed us to apply a lot of brain power in parallel, and when one branch of the project would slow down, another might speed up. Finally, we had an extremely capable project leader, one who could see the whole picture and was very adept at organizing and motivating people.

I don't mean to suggest that these attributes are all necessary ingredients for a successful Polymath project, but I think it is fair to say that, at least in this case, they were sufficient.


\section{Wouter Castryck}

\emph{Wouter Castryck is a postdoctoral fellow of FWO-Vlaanderen.}

In June last year, one of my colleagues informed me about the Polymath8 project and, in particular, about the programming challenge of finding admissible $k$-tuples whose diameter $H$ is as small as possible. I decided to give it a try and join ``team $H$'' of the production line. Luckily, I jumped in shortly after the project had started, at a point where there was still some low-hanging fruit. Like others I experienced how addictive it was to search for smaller values of $H$, while trying new computational tricks and combining them with ideas of the other participants. Because the value of $k$ kept decreasing as well, new challenges popped up every other day or so, which fueled the excitement. The whole event was intense and highly interactive, and progress was made at an incredible speed. (It is not academic to say so, but when the other teams managed to decrease $k$ to a size where our computational methods became superfluous, this was a bit of a disappointment. Luckily, Maynard's work on prime triples, quadruples, etc. put larger values of $k$ back into play.)

Along the way it became clear that the online arena in which it all took place attracted many spectators, and it felt like a privilege to be in there. It did require a mental switch to post naive (and sometimes wrong) ideas on the public forum, but in the end I agree with Andrew Sutherland that this is not necessarily a bad thing.

Gradually, by reading the blog posts, I also learned about the other parts of Zhang's proof. On a personal level, this I found the most enriching bit: participating in the Polymath8 project was a very stimulating way of learning and appreciating a new part of mathematics. This is definitely thanks to the clarifying and enthusing way in which Terence Tao administered the project. At the same time, I must admit that I did not grasp every detail from A to Z. From the point of view of a ``coauthor'' this is somewhat uncomfortable, but it may be inherent to the production line model along which the polymath projects are organized.


\section{Emmanuel Kowalski}

\emph{Emmanuel Kowalski is a professor at ETH.}

I remember reading the first blog posts of T. Gowers concerning the
Polymath idea, and following briefly the first steps at that time.  I
was intrigued but did not participate. As far as I remember, my only
public comment was a brief note after that progress was successful,
which still illustrates how I feel about the way this idea turned out:
it provides a striking, largely publicly available, illustration of
the way mathematical research works (or at least, often works...)
Although Gowers envisioned a collaboration with a very large number of
participants, what actually happened felt in fact closer to a
``standard'' collaboration, with the important difference that anybody
was free to jump in at any time, and to become a collaborator.
\par
I also spent a little time thinking about a later Polymath project
(finding large primes deterministically), but this also did not reach
the level of actual participation. Thus, the Polymath8 project is the
first time I have actually joined actively in this type of
mathematical work.
\par
In a general way, I feel that the Polymath idea does not fit very well
the way I tend to work -- this is not a criticism, but a statement of
fact.  In particular, following closely a free-for-all discussion on a
blog requires a type of focus and concentration which I often find
difficult to keep over more than a few days.  
\par
In the case of the Polymath8 project, however, a few coincidences led
to my participation:
\par
\begin{itemize}
\item I knew much of the subject material quite well, having studied
  the Goldston-Pintz-Y\i ld\i r\i m method relatively closely, as well
  as understanding (with varying degrees of depth) many of the
  ingredients of Zhang's work;
\item I had the occasion, during conferences and other visits, to
  discuss Zhang's work intensively with some of the best experts in
  the first weeks following the appearance of his paper; I reported on
  some of these discussions on my own blog, but I did not, at that
  time, track the Polymath8 project closely;
\item In particular, with \'E. Fouvry, Ph. Michel and P. Nelson, we
  were quickly able to understand, rework, and finally improve the
  crucial ``Type III'' estimates that were closely related to the work
  of Friedlander and Iwaniec concerning the ternary divisor function,
  which Fouvry, Michel and I had independently improved quite strongly
  a few months before; here the main point was the use of Deligne's
  formalism of the Riemann Hypothesis over finite fields, which we
  understood much better through the theory of trace functions over
  finite fields;
\item The final point was Ph. Michel's visit to California in the
  Summer 2013, during which it was formally decided to develop and
  incorporate this result in the Polymath8 project.
\end{itemize}

In the end, my closest involvement with the project came when the
first draft of the paper was written-up (I would probably have written
the section on Deligne's formalism, if Ph. Michel had not begun it
first) and time came to review and check it.  Since I was formally
identified as a participant to the project, I felt strongly the need
to read through the whole text, especially since parts of the claims
were really dependent on numerical facts for their correctness and
interest.  At first, I was only thinking that I would do a quick
review, but I find it difficult to read papers this way, and I ended
up doing a line-by-line check of the whole text, checking all
computations and reproducing (without copy-paste) all numerical
constants involved. It turned out that I only found one rather minor
mathematical mistake in the original text, and in particular that no
numerical value needed to be changed.  I still find this rather
remarkable.  Interestingly, the referee reports, which were extremely
thorough, turned out a few more such slips.  This illustrates how hard
it is to get a long mathematical argument rigorously right!
\par
Because this review was also done largely in public (on T. Tao's blog,
where the referee reports were also made available by public links) or
semi-publicly (in the Dropbox folder where the paper was written, and
in my own internal version-control setup), I think this can provide
interesting insights concerning the way mistakes may arise and be
corrected in long mathematical papers.  I feel however that the number
of problems was exceptionally small. In fact, I believe that a
similarly long paper of Polymath style, written without a strong
guiding hand like Tao's, could face serious issues or questions
concerning its correctness, if it turned out that sections were
written very independently, with authors who are working in a loose
organization, and do not fully master every part of the arguments.
This seems to me one important issue to keep in mind if this type of
collaborative work becomes more commonplace.
\par
Finally, I would like to emphasize a very positive point that was made
by R. de la Bret\'eche in his report \cite{breteche} on the Polymath8a project: this project unified
what might have been otherwise many small-scale improvements and
tweaks to Zhang's work, which would certainly have appeared much more
slowly -- if at all, in view of Maynard's breakthrough!


\section{Philippe Michel}

\emph{Philippe Michel is a Professor at EPF Lausanne.}

My involvement to the Polymath8 project began as a passive blog reader; when Zhang's breakthrough paper was publicly released, we, together with Paul Nelson in Lausanne, just like many other mathematicians, started to read the paper avidly. I should say that the material was not entirely a stranger to me --my master thesis subject was to understand the Acta Mathematica paper of Bombieri-Friedlander-Iwaniec and works related to it like the Inventiones paper of Deshouillers-Iwaniec-- but during the previous years, my mathematical interests had been rather remote and it was only during the past few months that I have been leaning back towards this area of analytic number theory - the study of the distribution of arithmetic functions along primes and arithmetic progressions - due to recent joint work with Etienne Fouvry and Emmanuel Kowalski. So I was then a bit rusty on Linnik's dispersion method and reading Terry Tao's blog posts on Zhang's work and the synthesis of the subsequent progress made by Polymath8 was tremendously useful to freshen up my memories and to test my intuitions. I should say that besides this, I found also very useful and instructive to read the comments made by other people at the end of each blog post.

My active participation to the project began only near the middle of June 2013; at that time, the Polymath8 project was already going full steam and had already made remarkable advances: the most visible one was a significant improvement over Zhang's initial constant $7 \times 10^7$ (dividing it by more than $1000$); another impressive result was the proof of an equidistribution theorem for primes in arithmetic progressions of large smooth moduli, going beyond the Bombieri--Vinogradov range, whose proof used ``only'' Weil's bound for exponential sums in one variable (Deligne's work being replaced by several applications of Cauchy inequality and of the $q$-van der Corput method). 

One portion of the argument which was still hurting the exponents concerned the ``type III sums''. We\footnote{E. Fouvry, E. Kowalski, P. Nelson, and myself.} had realized at an early stage that our previous work on the distribution of the ternary divisor function in large arithmetic progressions might help the cause; however it was only during the Fouvry 60 conference at CIRM in the middle of June, that we made this concrete and decided it could be worth contacting Polymath (through Terry Tao). It is highly probable that without the very existence of the Polymath8 project, its openness and the possibility for anybody to bring its own contribution, we would never have dared to make such a technical improvement public and, most likely, we would have forgotten about it after some time! This work was also an occasion to develop material of broader interest: one can find there a general account on $\ell$-adic trace functions which hopefully will be useful to the working analytic number theorist; for instance, we provide an easy to use presentation of the $q$-van der Corput method for general trace functions.  Another interesting outcome is the following: the quest for improvements on the numerical value of the distribution exponent has led to quite sophisticated transformations of the sums appearing in the dispersion method. The fact that the resulting ``complete'' algebraic exponential sums and their associated sheaves have a nice geometrical structure is a pleasant surprise which triggers a lot of interesting questions on the $\ell$-adic side.

I like very much working in a collaborative manner\footnote{The last non-survey paper I have written alone goes back to 2004.} possibly with fairly large groups of people; yet, this first participation to a Polymath project was at an entirely different scale and it took me some time to adapt.  One issue was to cope with the continuous flow of comments and new ideas made by the many participants (in particular, I often wondered whether I was contributing enough by comparison with others); another was to absorb the fact that the project was performed under the public eyes (and the vertiginous feeling that any mistake could be known to any mathematician in the world and would stay forever). Like for others, these concerns eventually disappeared and I began to fully enjoy the spontaneity of having everybody working openly through a public blog; at some point this even became slightly addictive: while attending a thesis defense (not in my area fortunately!), I found myself checking the latest development of the project and sending an email to some of the Polymath people to test an idea I just had on how to handle some unpleasant exponential sum. All in all, this has been a fun and rewarding experience and I am very thankful to Terry Tao for setting up and conducting the project; I am also thankful to the other participants for their constructive and highly collaborative attitude.

\end{document}